\documentclass[11pt]{article}
\usepackage{amssymb,amsfonts,amsmath, psfrag}
\usepackage{graphicx, colordvi, cite}
\usepackage{mathrsfs}
\usepackage{rotating}
\usepackage{hyperref}
\usepackage{lmodern}
\usepackage{bm}
\usepackage{amsthm}
\usepackage{url}
\newtheorem{theo}{Theorem}

\makeatletter \@addtoreset{equation}{section}
\@addtoreset{theo}{section} \makeatother

\usepackage{geometry}

\usepackage{graphicx}
\usepackage{tikz}
\usetikzlibrary{matrix,arrows}
\usetikzlibrary{positioning}
\usetikzlibrary{fit}
\usetikzlibrary{patterns}

\newcommand{\pattern}[4]{										% mesh pattern
	\raisebox{0.6ex}{
		\begin{tikzpicture}[scale=0.35, baseline=(current bounding box.center), #1]
		\foreach \x/\y in {#4}		\fill[gray!20] (\x,\y) rectangle +(1,1);
		\draw (0.01,0.01) grid (#2+0.99,#2+0.99);
		\foreach \x/\y in {#3}		\filldraw (\x,\y) circle (6pt);
		\end{tikzpicture}}
}

\definecolor{red}{rgb}{1,0,0}
{}
\begin{document}
\begin{center}
{\bf \large The last distribution-equivalence class of mesh patterns of length 2}
\vspace{0.3cm}

\end{center}
\begin{center}

 Zuo-Ru Zhang, \ Hongkuan Zhao\\[8pt]
School of Mathematical Science\\
Hebei Normal University\\
 Shijiazhuang 050024, P. R. China\\[6pt]

{\tt  zrzhang@hebtu.edu.cn}

{\tt  hkzhao2001@163.com}

\end{center}

\vspace{0.3cm} \noindent{\bf Abstract.} We prove the remaining conjectural distribution-equivalence of mesh patterns of length $2$ from 2019. This completes the classification of distribution-equivalence and Wilf-equivalence classes for mesh patterns of length $2$. Consequently the numbers of such classes are $105$ and $46$, respectively. We approach this by constructing a marked occurrence bijection, in the spirit of the cluster method for consecutive patterns.

\vspace{0.3cm}
\noindent{\bf Keywords: } mesh pattern $\cdot$ marked occurrence $\cdot$ equidistribution $\cdot$ Wilf-equivalence

\noindent {\bf AMS Subject Classifications:}  05A05, 05A15.
\section{Introduction}
%Denote by $[a,b]:=\{a,a+1,\ldots,b\}$ and $\mathcal{S}_n$ the set of permutations of $[n]:=[1,n]$. A (classical permutation) pattern is a permutation. Given $\sigma\in \mathcal{S}_n$ and $\pi\in \mathcal{S}_m$, a subsequence of the entries of $\sigma$ is called an occurrence of the pattern $\pi$ if it has the same relative order as all of the entries of $\pi$.
Denote by $[a,b]$ the interval of the integers from $a$ to $b$. Recall that a mesh pattern of length $k$ is a pair $(\tau,R)$, where $\tau$ is a permutation of length $k$ and $R$ is a subset of $[0,k]\times[0,k]$. Let $(i,j)$ denote the box whose corners have coordinates $(i,j)$, $(i,j+1)$, $(i+1,j+1)$ and $(i+1,j)$.
Mesh patterns can be drawn by shading the boxes in $R$. The picture 
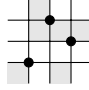
\[
\pattern{scale=0.8}{3}{1/1,2/3,3/2}{0/0,1/2,1/3,2/0,2/2,3/2}
\]
represents the mesh pattern with $\tau=132$ and $R=\{(0,0),(1,2),(1,3),(2,0),(2,2),(3,2)\}$. 

A subsequence $\pi'$ of a permutation $\pi$ is an occurrence of a mesh pattern $(\tau,R)$ if ($1$) $\pi'$ is an occurrence of $\tau$, and ($2$) the shaded squares specified by $R$ contain no elements of $\pi$ other than those in $\pi'$.
Throughout this paper, when we say that $ab$ is an occurrence of a mesh pattern of length $2$, $a$ and $b$ refer to the values of the underlying (classical) pattern.
Let $\operatorname{Occ}_p(\pi)$ denote the set of all occurrences of the pattern $p$ in $\pi$ and $o_p(\pi)$ the cardinality of $\operatorname{Occ}_p(\pi)$. Two patterns $p$ and $q$ are equidistributed if, for all $n\geq 0$, $F_n^p(y):=\sum_{\pi\in\mathcal{S}_n}y^{o_p(\pi)}=F_n^q(y):=\sum_{\pi\in\mathcal{S}_n}y^{o_q(\pi)}$. In particular, $p$ and $q$ are Wilf-equivalent if $F_n^p(0)=F_n^q(0)$. The corresponding equivalence classes are called distribution-equivalence classes and Wilf-equivalence classes, respectively.

The notion of mesh patterns was introduced by Br\"and\'en and Claesson~\cite{BC} to study the explicit expansions for certain permutation statistics as linear combinations of classical permutation patterns. Hilmarsson et al.~\cite{Hil} started the systematic study of avoidance of mesh patterns of length $2$. Kitaev and Zhang~\cite{KZ} further studied the distributions of mesh patterns.  Significant progress on determining the distribution-equivalence and Wilf-equivalence classes for mesh patterns of length $2$ has been made by many authors, see~\cite{FFKLSS,HZ,Hil,KZ,KZZ,SKZ}. To date, the only remaining conjectural equidistributed class is Class $69$ comprising of the patterns in the set $\{\pattern{scale=0.5}{2}{1/1,2/2}{1/2,1/1,2/1,0/0},\pattern{scale=0.5}{2}{1/1,2/2}{2/2,0/1,1/1,1/0},\pattern{scale=0.5}{2}{1/1,2/2}{0/2,1/1,2/1,1/0},\pattern{scale=0.5}{2}{1/1,2/2}{1/2,0/1,1/1,2/0} \}$, where the first two patterns, as well as the last two patterns, are trivially equivalent~\cite[Remark~1.3]{SKZ}.

The goal of this paper is to prove this conjecture.
Equivalently, we prove the following theorem.
\begin{theo}\label{57-58-equivalent}
The patterns $P=\pattern{scale = 0.6}{2}{1/1,2/2}{0/1,1/1,1/2,2/0}$ and $Q=\pattern{scale = 0.6}{2}{1/2,2/1}{0/1,1/1,1/2,2/0}$ are equidistributed.
\end{theo}

Together with the previously known equivalences and separations, this completes the classification of both distribution-equivalence and Wilf-equivalence classes for mesh patterns of length $2$. Consequently, there are exactly $105$ distribution-equivalence classes and $46$ Wilf-equivalence classes.

Our proof is bijective.
Instead of constructing a direct occurrence-preserving bijection on permutations, we use marked occurrences, inspired by the cluster method for consecutive patterns~\cite{EliNoy}.
More precisely, let
\[
\mathcal{MOS}_n(p):=\{(\pi;M)\mid \pi \in \mathcal{S}_n, M\subseteq \operatorname{Occ}_p(\pi)\}.
\]

Then 
\begin{align*}
F_n^p(y)
 &=\sum_{\pi\in\mathcal{S}_n}y^{o_p(\pi)} \\
 &=\sum_{\pi\in\mathcal{S}_n}(1+z)^{o_p(\pi)} \qquad (y=1+z) \\
 &=\sum_{\pi\in\mathcal{S}_n}\sum_{k=0}^{o_p(\pi)}\binom{o_p(\pi)}{k}z^k\\
 &=\sum_{\pi\in\mathcal{S}_n}\sum_{M\subseteq \operatorname{Occ}_p(\pi)}z^{|M|}\\
 &=\sum_{(\pi;M) \in \mathcal{MOS}_n(p)}z^{|M|}.
\end{align*}
For $(\pi;M)\in\mathcal{MOS}_n(p)$, we call $M$ a marked occurrence set of $p$ in $\pi$.
Therefore, $F_n^P(y)=F_n^Q(y)$ if there is a bijection $\phi: \mathcal{MOS}_n(P)\rightarrow\mathcal{MOS}_n(Q)$ that preserves the cardinality of the marked occurrence set.
The rest of the paper is devoted to the construction of this bijection.
\section{The construction of \texorpdfstring{$\phi$}{phi}}
\textbf{Observation:} Suppose $ab$ and $cd$ are two distinct occurrences of $P$ in $\pi$, then none of the following cases happens.

\noindent
\(
\begin{aligned}[t]
\text{(i)}\quad & a=c;\\
\text{(ii)}\quad & b=d;\\
\text{(iii)}\quad &
\pi^{-1}(a)<\pi^{-1}(c)<\pi^{-1}(b)<\pi^{-1}(d)
\quad\text{or}\quad
\pi^{-1}(c)<\pi^{-1}(a)<\pi^{-1}(d)<\pi^{-1}(b).
\end{aligned}
\)

The proof is straightforward and we omit it.

Now we construct $\phi$.

Let $(\pi;M)\in\mathcal{MOS}_n(P)$. We call a subset of $M$ of the form $\{a_1a_2,a_2a_3,\ldots,a_{m-1}a_m\}$ a $P$-path, where $a_1a_2$ and $a_{m-1}a_m$ are, respectively, the unique occurrences of $P$ that contain $a_1$ and $a_m$ in $M$.

(1) We first assume that $M$ is just a $P$-path, namely,
\[
M=\{a_1a_2,a_2a_3,\ldots,a_{m-1}a_m\}.
\]
Clearly $a_1<a_2<\ldots<a_m$.

Modify $\pi$ by simultaneously performing the following substitutions:
\[
a_1 \to a_m,\, a_2 \to a_1,\, \ldots,\, a_{m-1} \to a_{m-2},\, a_m \to a_{m-1}.
\]
Then we obtain $(\sigma;N)$ where $N=\{a_ma_1,a_ma_2,\ldots,a_ma_{m-1}\}$. To show $(\sigma;N)\in\mathcal{MOS}_n(Q)$, i.e., $N\subseteq \operatorname{Occ}_{Q}(\sigma)$, we just need to check the following three conditions.

\textbf{Left} \,Let $x>a_1$ be an entry to the left of $a_1$ in $\pi$.
Since $a_1a_2$ is an occurrence of the pattern $P$, the defining shading condition implies that $x>a_2$.
Since $a_2a_3$ is also an occurrence of the pattern $P$, the same argument gives $x>a_3$.
Iterating this argument, we obtain $x>a_m$.

\textbf{Middle} \,In $\pi$, the entries between $a_i$ and $a_{i+1}$ ($1\leq i \leq m-1$) are less than $a_i$.
Thus the entries between $a_1$ and $a_{j+1}$ (except for $a_j$, $1\leq j \leq m-1$) are less than $a_{j}$.
In $\sigma$, the entries between $a_m$ and $a_j$ are exactly the entries between $a_1$ and $a_{j+1}$ in $\pi$, with $a_j$ replaced by $a_1$.
So they are less than $a_j$ (and $a_m$).

\textbf{Right} \,In $\sigma$, the entries to the right of $a_k$ are exactly those to the right of $a_{k+1}$ in $\pi$, with $a_m$ replaced by $a_{k+1}$ ($1\leq k\leq m-1$).
Since $a_ka_{k+1}$ is an occurrence of the pattern $P$, the entries to the right of $a_{k+1}$ are greater than $a_k$.
Thus in $\sigma$, the entries to the right of $a_k$ are greater than $a_k$.

Therefore, $N\subseteq \operatorname{Occ}_{Q}(\sigma)$.

(2) Now assume that $b_1b_2,b_2b_3,\ldots,b_{s-1}b_s$ and $c_1c_2,\ldots,c_{t-1}c_t$ are two disjoint $P$-paths (otherwise they form a new $P$-path).
We show that the replacements associated with them don't interfere with one another.

If $\pi^{-1}(b_s)<\pi^{-1}(c_1)$, the argument is clear.

If $\pi^{-1}(b_i)<\pi^{-1}(c_1)<\pi^{-1}(c_t)<\pi^{-1}(b_{i+1})$ for some $i$, then treating the path $c$ clearly doesn't affect the path $b$.
Since $b_ib_{i+1}$ is an occurrence of the pattern $P$, $c_1<c_t<b_i<b_{i+1}\leq b_s$. When treating the path $b$, we merely put $b_i$ to the right of $c_t$ while putting $b_s$ to the left of $c_1$; this doesn't affect the path $c$.

Therefore, for a general $(\pi;M)\in\mathcal{MOS}_n(P)$, we treat every $P$-path like above and then obtain $\phi[(\pi;M)]$.
The inverse map is obtained analogously, so $\phi$ is the desired bijection; we omit the details here.
%Analogously one can verify that $\phi$ is invertible and thus $\phi$ is the bijection we desired, we omit the details here.

We use the figure below to illustrate the map $\phi$.

\begin{figure}[htbp]
\centering

\[
\resizebox{\linewidth}{!}{%
\begin{tikzpicture}[
    every node/.style={inner sep=1pt, font=\large},
    dot/.style={font=\large},
    arc/.style={line width=0.55pt}
]
% baseline
\def\y{0}

% nodes
\node (a1) at (0,\y) {$a_1$};
\node (a2) at (2,\y) {$a_2$};
\node (b1) at (4,\y) {$b_1$};
\node (b2) at (5.8,\y) {$b_2$};
\node (a3) at (8,\y) {$a_3$};
\node (a4) at (10,\y) {$a_4$};
\node (c1) at (12.2,\y) {$c_1$};
\node (c2) at (14.2,\y) {$c_2$};
\node (c3) at (16.2,\y) {$c_3$};

% dots
\node[dot] at (-0.9,\y) {$\cdots$};
\node[dot] at (1,\y) {$\cdots$};
\node[dot] at (3,\y) {$\cdots$};
\node[dot] at (4.9,\y) {$\cdots$};
\node[dot] at (6.9,\y) {$\cdots$};
\node[dot] at (9,\y) {$\cdots$};
\node[dot] at (11.1,\y) {$\cdots$};
\node[dot] at (13.2,\y) {$\cdots$};
\node[dot] at (15.2,\y) {$\cdots$};
\node[dot] at (17.3,\y) {$\cdots$};

% arcs
\draw[arc] (a1.south) .. controls (0.35,-1.00) and (1.65,-1.00) .. (a2.south);
\draw[arc] (b1.south) .. controls (4.25,-0.75) and (5.55,-0.75) .. (b2.south);
\draw[arc] (a3.south) .. controls (8.35,-0.95) and (9.65,-0.95) .. (a4.south);
\draw[arc] (c1.south) .. controls (12.55,-0.95) and (13.85,-0.95) .. (c2.south);
\draw[arc] (c2.south) .. controls (14.55,-0.95) and (15.85,-0.95) .. (c3.south);

\draw[arc] (a2.south) .. controls (2.9,-1.95) and (7.1,-1.95) .. (a3.south);
\end{tikzpicture}%
}
\]
\vspace{1em}
\[
\resizebox{\linewidth}{!}{%
\begin{tikzpicture}[
    every node/.style={inner sep=1pt, font=\large},
    dot/.style={font=\large},
    arc/.style={line width=0.55pt}
]
% baseline
\def\y{0}

% nodes
\node (a4) at (0,\y) {$a_4$};
\node (a1) at (2,\y) {$a_1$};
\node (b2) at (4.3,\y) {$b_2$};
\node (b1) at (6.1,\y) {$b_1$};
\node (a2) at (8.4,\y) {$a_2$};
\node (a3) at (10.6,\y) {$a_3$};
\node (c3) at (13.1,\y) {$c_3$};
\node (c1) at (15.2,\y) {$c_1$};
\node (c2) at (17.3,\y) {$c_2$};

% dots
\node[dot] at (-0.9,\y) {$\cdots$};
\node[dot] at (1,\y) {$\cdots$};
\node[dot] at (3.15,\y) {$\cdots$};
\node[dot] at (5.2,\y) {$\cdots$};
\node[dot] at (7.25,\y) {$\cdots$};
\node[dot] at (9.5,\y) {$\cdots$};
\node[dot] at (11.8,\y) {$\cdots$};
\node[dot] at (14.15,\y) {$\cdots$};
\node[dot] at (16.25,\y) {$\cdots$};
\node[dot] at (18.35,\y) {$\cdots$};

% arcs
\draw[arc] (a4.south) .. controls (0.35,-0.95) and (1.65,-0.95) .. (a1.south);
\draw[arc] (b2.south) .. controls (4.55,-0.75) and (5.85,-0.75) .. (b1.south);
\draw[arc] (c3.south) .. controls (13.45,-0.95) and (14.85,-0.95) .. (c1.south);

\draw[arc] (a4.south) .. controls (1.25,-2.15) and (7.15,-2.15) .. (a2.south);
\draw[arc] (a4.south) .. controls (1.45,-3.00) and (9.15,-3.00) .. (a3.south);
\draw[arc] (c3.south) .. controls (13.8,-2.15) and (16.6,-2.15) .. (c2.south);
\end{tikzpicture}%
}
\]
\caption{$(\pi;\{a_1a_2,a_2a_3,a_3a_4,b_1b_2,c_1c_2,c_2c_3\})$ goes to $(\sigma;\{a_4a_1,a_4a_2,a_4a_3,b_2b_1,c_3c_1,c_3c_2\})$ under $\phi$.}
\end{figure}
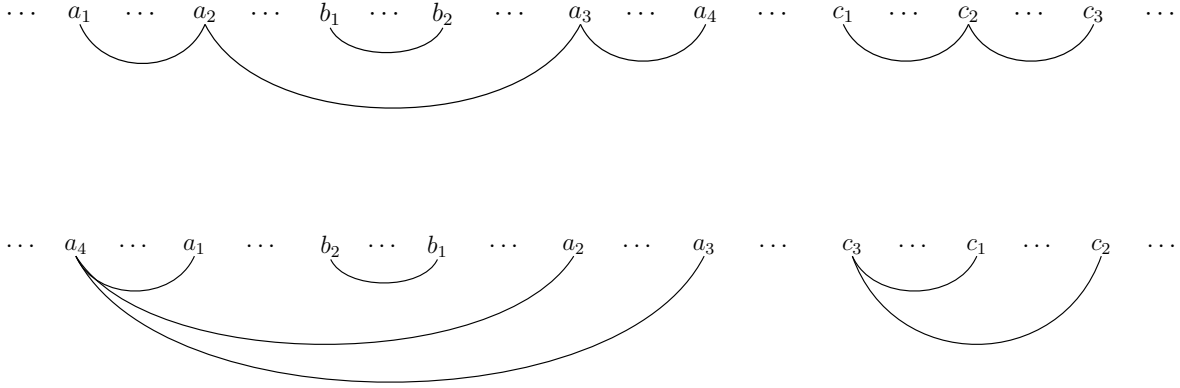

\end{document}